\documentclass[11pt]{amsart}    
\usepackage{latexsym, amsmath, amsthm, amssymb, setspace, verbatim}
\usepackage{longtable}
\usepackage[applemac]{inputenc}
\usepackage{hyperref}

\title{ The Rokhlin dimension of topological $\IZ^m$-actions }
\author{ Gábor Szabó }
\address{Westfälische Wilhelms-Universität, Fachbereich Mathematik, \phantom{--------------}\linebreak \text{}\hspace{3.5mm} Einsteinstrasse 62, 48149 Münster, Germany}
\email{gabor.szabo@uni-muenster.de}
\thanks{\emph{Supported by:} SFB 878 \emph{Groups, Geometry and Actions} and GIF Grant 1137-30.6/2011}

\begin{document}

\renewcommand\matrix[1]{\left(\begin{array}{*{10}{c}} #1 \end{array}\right)}  
\newcommand\set[1]{\left\{#1\right\}}  
\newcommand\mset[1]{\left\{\!\!\left\{#1\right\}\!\!\right\}}

\newcommand{\IA}[0]{\mathbb{A}} \newcommand{\IB}[0]{\mathbb{B}}
\newcommand{\IC}[0]{\mathbb{C}} \newcommand{\ID}[0]{\mathbb{D}}
\newcommand{\IE}[0]{\mathbb{E}} \newcommand{\IF}[0]{\mathbb{F}}
\newcommand{\IG}[0]{\mathbb{G}} \newcommand{\IH}[0]{\mathbb{H}}
\newcommand{\II}[0]{\mathbb{I}} \renewcommand{\IJ}[0]{\mathbb{J}}
\newcommand{\IK}[0]{\mathbb{K}} \newcommand{\IL}[0]{\mathbb{L}}
\newcommand{\IM}[0]{\mathbb{M}} \newcommand{\IN}[0]{\mathbb{N}}
\newcommand{\IO}[0]{\mathbb{O}} \newcommand{\IP}[0]{\mathbb{P}}
\newcommand{\IQ}[0]{\mathbb{Q}} \newcommand{\IR}[0]{\mathbb{R}}
\newcommand{\IS}[0]{\mathbb{S}} \newcommand{\IT}[0]{\mathbb{T}}
\newcommand{\IU}[0]{\mathbb{U}} \newcommand{\IV}[0]{\mathbb{V}}
\newcommand{\IW}[0]{\mathbb{W}} \newcommand{\IX}[0]{\mathbb{X}}
\newcommand{\IY}[0]{\mathbb{Y}} \newcommand{\IZ}[0]{\mathbb{Z}}

\newcommand{\CA}[0]{\mathcal{A}} \newcommand{\CB}[0]{\mathcal{B}}
\newcommand{\CC}[0]{\mathcal{C}} \newcommand{\CD}[0]{\mathcal{D}}
\newcommand{\CE}[0]{\mathcal{E}} \newcommand{\CF}[0]{\mathcal{F}}
\newcommand{\CG}[0]{\mathcal{G}} \newcommand{\CH}[0]{\mathcal{H}}
\newcommand{\CI}[0]{\mathcal{I}} \newcommand{\CJ}[0]{\mathcal{J}}
\newcommand{\CK}[0]{\mathcal{K}} \newcommand{\CL}[0]{\mathcal{L}}
\newcommand{\CM}[0]{\mathcal{M}} \newcommand{\CN}[0]{\mathcal{N}}
\newcommand{\CO}[0]{\mathcal{O}} \newcommand{\CP}[0]{\mathcal{P}}
\newcommand{\CQ}[0]{\mathcal{Q}} \newcommand{\CR}[0]{\mathcal{R}}
\newcommand{\CS}[0]{\mathcal{S}} \newcommand{\CT}[0]{\mathcal{T}}
\newcommand{\CU}[0]{\mathcal{U}} \newcommand{\CV}[0]{\mathcal{V}}
\newcommand{\CW}[0]{\mathcal{W}} \newcommand{\CX}[0]{\mathcal{X}}
\newcommand{\CY}[0]{\mathcal{Y}} \newcommand{\CZ}[0]{\mathcal{Z}}

\newcommand{\FA}[0]{\mathfrak{A}} \newcommand{\FB}[0]{\mathfrak{B}}
\newcommand{\FC}[0]{\mathfrak{C}} \newcommand{\FD}[0]{\mathfrak{D}}
\newcommand{\FE}[0]{\mathfrak{E}} \newcommand{\FF}[0]{\mathfrak{F}}
\newcommand{\FG}[0]{\mathfrak{G}} \newcommand{\FH}[0]{\mathfrak{H}}
\newcommand{\FI}[0]{\mathfrak{I}} \newcommand{\FJ}[0]{\mathfrak{J}}
\newcommand{\FK}[0]{\mathfrak{K}} \newcommand{\FL}[0]{\mathfrak{L}}
\newcommand{\FM}[0]{\mathfrak{M}} \newcommand{\FN}[0]{\mathfrak{N}}
\newcommand{\FO}[0]{\mathfrak{O}} \newcommand{\FP}[0]{\mathfrak{P}}
\newcommand{\FQ}[0]{\mathfrak{Q}} \newcommand{\FR}[0]{\mathfrak{R}}
\newcommand{\FS}[0]{\mathfrak{S}} \newcommand{\FT}[0]{\mathfrak{T}}
\newcommand{\FU}[0]{\mathfrak{U}} \newcommand{\FV}[0]{\mathfrak{V}}
\newcommand{\FW}[0]{\mathfrak{W}} \newcommand{\FX}[0]{\mathfrak{X}}
\newcommand{\FY}[0]{\mathfrak{Y}} \newcommand{\FZ}[0]{\mathfrak{Z}}

\newcommand{\Ra}[0]{\Rightarrow}
\newcommand{\La}[0]{\Leftarrow}
\newcommand{\LRa}[0]{\Leftrightarrow}

\renewcommand{\phi}[0]{\varphi}
\newcommand{\eps}[0]{\varepsilon}

\newcommand{\quer}[0]{\overline}
\newcommand{\uber}[0]{\choose}
\newcommand{\ord}[0]{\operatorname{ord}}		
\newcommand{\GL}[0]{\operatorname{GL}}
\newcommand{\supp}[0]{\operatorname{supp}}	
\newcommand{\id}[0]{\operatorname{id}}		
\newcommand{\Sp}[0]{\operatorname{Sp}}		
\newcommand{\eins}[0]{\mathbf{1}}			
\newcommand{\diag}[0]{\operatorname{diag}}
\newcommand{\stark}[0]{\stackrel{s}{\to}}
\newcommand{\schwach}[0]{\stackrel{w}{\to}}
\newcommand{\ind}[0]{\operatorname{ind}}
\newcommand{\auf}[1]{\quad\stackrel{#1}{\longrightarrow}\quad}
\newcommand{\hull}[0]{\operatorname{hull}}
\newcommand{\prim}[0]{\operatorname{Prim}}
\newcommand{\ad}[0]{\operatorname{Ad}}
\newcommand{\quot}[0]{\operatorname{Quot}}
\newcommand{\ext}[0]{\operatorname{Ext}}
\newcommand{\ev}[0]{\operatorname{ev}}
\newcommand{\fin}[0]{{\subset\!\!\!\subset}}
\newcommand{\diam}[0]{\operatorname{diam}}
\newcommand{\Hom}[0]{\operatorname{Hom}}
\newcommand{\Aut}[0]{\operatorname{Aut}}
\newcommand{\del}[0]{\partial}
\newcommand{\dimnuc}[0]{\dim_{\mathrm{nuc}}}
\newcommand{\dimrok}[0]{\dim_{\mathrm{Rok}}}
\newcommand{\dimrokcyc}[0]{\dim_{\mathrm{Rok}}^{\mathrm{cyc}}}
\newcommand{\dimrokcycc}[0]{\dim_{\mathrm{Rok}}^{\mathrm{cyc,c}}}
\newcommand{\dr}[0]{\operatorname{dr}}
\newcommand*\onto{\ensuremath{\joinrel\relbar\joinrel\twoheadrightarrow}} 
\newcommand*\into{\ensuremath{\lhook\joinrel\relbar\joinrel\rightarrow}}  
\newcommand{\cstar}[0]{$\mathrm{C}^*$}

\newtheorem{theorem}{Theorem}[section]
\newtheorem{cor}[theorem]{Corollary}
\newtheorem{lemma}[theorem]{Lemma}
\newtheorem{prop}[theorem]{Proposition}
\newtheorem*{theoreme}{Theorem}

\theoremstyle{definition}
\newtheorem{defi}[theorem]{Definition}
\newtheorem{nota}[theorem]{Notation}
\newtheorem{rem}[theorem]{Remark}
\newtheorem*{reme}{Remark}
\newtheorem{question}[theorem]{Question}

\begin{abstract}
We study the topological variant of Rokhlin dimension for topological dynamical systems $(X,\alpha,\IZ^m)$ in the case where $X$ is assumed to have finite covering dimension. Finite Rokhlin dimension in this sense is a property that implies finite Rokhlin dimension of the induced action on \cstar-algebraic level, as was discussed in a recent paper by Hirshberg, Winter and Zacharias. In particular, it implies under these conditions that the transformation group \cstar-algebra has finite nuclear dimension. Generalizing partial results of Lindenstrauss and Gutman, we show that free $\IZ^m$-actions on finite dimensional spaces satisfy a strengthened version of the so-called marker property, which yields finite Rokhlin dimension for such actions.
\end{abstract}

\maketitle


\setcounter{section}{-1}
\section{Introduction}
\noindent
The study of group actions and their associated \cstar-algebras has always been a central theme in the theory of operator algebras. Topological dynamical systems in particular play a significant role in \cstar-algebra theory. 
The \cstar-algebras arising from topological dynamical systems, say from $(X,\alpha,G)$ for a countable discrete group $G$, a compact metric space $X$ and a continuous action $\alpha: G\curvearrowright X$, give rise to plenty of interesting questions about their structure. From the point of view of \cstar-algebra classification theory, this setting begs the question of whether the transformation group \cstar-algebras can be classified via their Elliott-invariant, at least under suitable conditions like simplicity.

By now there exist classification results for large classes of crossed products by $\IZ$-actions: early results of Putnam about characterizing crossed products of minimal homeomorphisms on the Cantor set as A$\IT$ algebras (see \cite{Putnam}) or of Elliott and Evans about irrational rotation algebras (see \cite{EllEv}) have set the stage for this project. In one of the more recent breakthroughs, Toms, Strung and Winter proved that crossed products of uniquely ergodic minimal homeomorphisms on infinite compact metrizable spaces with finite covering dimension are classified by ordered K-theory (see \cite{WiTo,StrWi}).

The case of $\IZ^m$-actions begins to gather more and more attention and interest as the next step for this long-term study of group actions. Towards a generalization of the celebrated results of \cite{GPS}, Giordano, Matui, Putnam and Skau undertook a deep study on topological orbit equivalence of Cantor minimal $\IZ^m$-actions in \cite{GMPS1} and \cite{GMPS2}. Going in a more similar direction as this paper, Phillips initiated the study of certain \emph{large subalgebras} inside crossed products by topological $\IZ^m$-actions as a means to study the radius of comparison of these \cstar-algebras (see \cite{Phi}). In particular, an important question concerning the \cstar-classification of these crossed products is when they have strict comparison of positive elements. 

The current trend in \cstar-algebra classification theory is to show certain regularity properties for a class of nuclear \cstar-algebras as the first intermediate step towards their $K$-theoretic classification. Following this philosophy, we focus on finite nuclear dimension (see \cite{WiZa}) as the regularity property of our choice and show that it is prevalent for transformation group \cstar-algebras of free $\IZ^m$-actions on finite dimensional spaces. Looking in the direction of the second step in this classification philosophy, Winter has found a method to deduce classification of simple transformation group \cstar-algebras of this type. This works under the prevalence of finite nuclear dimension and a technical condition involving the set of invariant ergodic measures (see \cite{Winter}). An alternative immediate application is to combine the main results of \cite{Lin} and \cite{MatuiSato} to obtain a classification result for free, minimal and uniquely ergodic $\IZ^m$-actions.

The concept of Rokhlin dimension developed in \cite{HWZ} is a natural tool for a systematic approach of showing finite nuclear dimension for crossed products. This notion has been introduced for finite group actions and integer actions on unital \cstar-algebras, but can be defined similarly for $\IZ^m$-actions. Integer actions with finite Rokhlin dimension have been shown to behave well with underlying \cstar-algebras of finite nuclear dimension. That is, the property of having finite nuclear dimension passes from the underlying \cstar-algebra to the crossed product. The aim of the first section is to define Rokhlin dimension for $\IZ^m$-actions on unital \cstar-algebras and establish the same permance property for $\IZ^m$-actions of finite Rokhlin dimension.

In the second section, we introduce a topological variant of Rokhlin dimension for topological $\IZ^m$-actions. In particular, finite Rokhlin dimension in the topological sense is designed in such a way that it implies finite Rokhlin dimension in the \cstar-algebraic sense for its induced \cstar-action. Moreover, it can be regarded as a topological analogue of the measure theoretic Rokhlin Lemma. The rest of the paper will be devoted to show that finite Rokhlin dimension in the topological sense is satisfied for free $\IZ^m$-action on finite dimensional spaces.

In the third section, we introduce a technical condition that is a stronger version of the small boundary property introduced by Lindenstrauss in \cite{Lind}, namely the bounded topological small boundary property. We show that free actions of countably infinite groups on finite dimensional spaces satisfy said property.

In the fourth section, we define the marker property as in \cite{Gutman} within the more general setting of countable group actions. We also define a stronger variant of this property for $\IZ^m$-actions, namely the controlled marker property. In \cite{Gutman}, Gutman has shown that aperiodic homeomorphisms on finite dimensional spaces have the marker property. We generalize this result and show that free actions of countably infinite groups have the marker property whenever the system has the bounded topological small boundary property. In the case of $\IZ^m$-actions, a closer look at the proof will even yield the controlled marker property under these conditions.

In the fifth section, we bring together the results of all the previous sections. Combining the third and fourth section, we can deduce that free $\IZ^m$-actions on finite dimensional spaces satisfy the controlled marker property. The consequence is finite Rokhlin dimension in the topological sense. By the second section, we get finite Rokhlin dimension for the induced \cstar-action. Applying the results of the first section, we get the main result:
{\theoreme Let $X$ be a compact metric space of finite covering dimension and let $\alpha: \IZ^m\curvearrowright X$ be a free continuous group action. Then the induced \cstar-action $\bar{\alpha}: \IZ^m\curvearrowright\CC(X)$ has finite Rokhlin dimension and the transformation group \cstar-algebra $\CC(X)\rtimes_{\bar{\alpha}}\IZ^m$ has finite nuclear dimension.}\vspace{3mm}

The results of this paper belong to the research of the author's PhD studies. 

\section{Preliminaries and \cstar-Rokhlin dimension}

\nota Throughout the whole paper, we will stick to the following notations unless specified otherwise:
\begin{itemize}
\item $X$ is a compact metric space. In applications, it is often assumed to have finite covering dimension.
\item $G$ is a countable and discrete group.
\item $A$ is a unital \cstar-algebra.
\item A completely positive map between \cstar-algebras is abbreviated \emph{c.p.} A completely positive and contractive map between \cstar-algebras is abbreviated \emph{c.p.c.}
\item Either $\alpha: G\curvearrowright X$ is a continuous group action or $\alpha: G\curvearrowright A$ is an action via automorphisms. When $\alpha: G\curvearrowright X$ is a topological action, we denote the induced \cstar-action by $\bar{\alpha}: G \curvearrowright\CC(X)$.  In applications, we often have $G=\IZ^m$.
\item Suppose $M$ is some set. If $F\subset M$ is a finite subset, we write $F\fin M$.
\item For $a,b\in A$ in some \cstar-algebra and $\eps>0$, $a=_\eps b$ means $\|a-b\|\leq\eps$.
\item Suppose we have $F,F'\fin A$ in a \cstar-algebra and $\eps>0$. The notation $F=_\eps F'$ means that for all $a\in F$, there exists $b\in F'$ with $a=_\eps b$ and vice versa.
\item Given $n\in\IN$, we define
\[B_n^m = \set{0,\dots,n-1}^m\subset\IZ^m.\]
Since $m$ will always be known from context to be the rank of $\IZ^m$, we will just write $B_n$ instead.
\item Similarly, we define
\[ J_n^m = J_n= \set{-n+1,\dots, 0, \dots , n}^m\subset\IZ^m.\]
\end{itemize}
We begin with a few remarks about almost order zero maps and nuclear dimension.

\defi Let $A,B$ be \cstar-algebras and $\delta\geq0$. A c.p.c.~map $\psi: A\to B$ has order zero up to $\delta$, if for all $a,b\in A$ with $ab=0$ we have $\|\psi(a)\psi(b)\|\leq\delta\|a\|\|b\|$. $\psi$ has order zero if $\delta=0$. \vspace{3mm}\\
Order zero maps in particular appear in the definition of nuclear dimension and decomposition rank (see \cite{KiWi,WiZa}),  which are important regularity properties for the classification theory of nuclear \cstar-algebras. We shall recall a stability property of order zero maps that is satisfied whenever the domain is finite dimensional (see \cite[2.5]{KiWi}).

\lemma \label{stability} Let $\CF$ be a finite dimensional \cstar-algebra. For all $\eps>0$, there exists $\delta>0$ with the following property:
Let $A$ be a \cstar-algebra and $\phi: \CF\to A$ a map of order zero up to $\delta$. Then there exists an order zero map $\phi': \CF\to A$ such that $\|\phi-\phi'\|\leq\eps$.

\defi[see \cite{WiZa}]
Let $A$ be a \cstar-algebra. $A$ is said to have nuclear dimension $n$, denoted by $\dimnuc(A)=n$, if $n$ is the smallest natural number with the following property:

For all $F\fin A$ and $\eps>0$, there exists a finite dimensional \cstar-algebra
\[\CF=\CF^{(0)}\oplus\dots\oplus\CF^{(n)}\]
and c.p. maps\quad $A\auf{\psi}\CF\auf{\phi} A$ \quad such that
\begin{itemize}
\item $\psi$ is a c.p.c.~map.
\item For all $i=0,\dots,n$, the map $\phi^{(i)} := \phi|_{\CF^{(i)}}$ is c.p.c.~order zero.
\item $\|\phi\circ\psi(a)-a\|\leq\eps$ for all $a\in F$.
\end{itemize}
In this context, the triple $(\CF,\psi,\phi)$ is called an $n$-decomposable c.p.~approximation of tolerance $\eps$ on $F$.
If no such $n$ exists, we write $\dimnuc(A)=\infty$.

\rem[see {\cite[Lemma A.4]{HWZ}} ] \label{orderdelta} Using \ref{stability}, one can show the following: 

If $n$ is a number such that for all $F\fin A$ and $\delta>0$, there exists a finite dimensional C*-algebra $\CF$ such that for all $\eta>0$, we can choose a c.p.~approximation $(\CF,\phi,\psi)$ of tolerance $\delta$ on $F$ with $\phi$ being decomposable into $n+1$ maps $\phi^{(i)}$ of order zero up to $\eta$, then $\dimnuc(A)\leq n$.

\defi \label{dimrok} 
Let $A$ be a unital \cstar-algebra, and let $\alpha: \IZ^m\curvearrowright A$ be a group action via automorphisms.
We say that the action $\alpha$ has (cyclic) Rokhlin dimension $d$, and write $\dimrokcyc(\alpha)=d$, if $d$ is the smallest natural number with the following property:

For all $F\fin A, \eps>0, n\in\IN$, there exist positive contractions 
$(f_v^{(l)})^{l=0,\dots,d}_{v\in B_n}$ in $A$ satisfying the following properties:
\begin{itemize}
\item[(1)] $\displaystyle \|\eins_A - \sum_{l=0}^d \sum_{v\in B_n} f_v^{(l)} \|\leq\eps$.
\item[(2)] $\|f_{v}^{(l)}f_{v'}^{(l)}\|\leq\eps$ for all $l=0,\dots,d$ and $v \neq v'$ in $B_n$.
\item[(3)] $\|\alpha^v(f_{w}^{(l)})-f_{v+w}^{(l)}\|\leq\eps$ for all $l=0,\dots,d$ and $v,w\in B_n$.
 \vspace{1mm}\\ (Note that $f_v^{(l)}$ denotes $f_{(v\mod n\IZ^m)}^{(l)}$, whenever $v\notin B_n$.)
\item[(4)] $\|[f^{(l)}_v, a]\|\leq\eps$ for all $l=0,\dots,d\; , v\in B_n$ and $a\in F$.
\end{itemize}
If there is no such $d$, we write $\dimrok(\alpha)=\infty$.\vspace{3mm}

\reme We define (cyclic) Rokhlin dimension with commuting towers of $\alpha$, written $\dimrokcycc(\alpha)$, in the same way with the additional property
\begin{itemize}
\item[(5)] $\|[f_v^{(l)},f_w^{(l')}]\|\leq\eps$ for all $l,l'=0,\dots,d$ and $v,w\in B_n$
\end{itemize}
for the Rokhlin elements $(f_v^{(l)})^{l=0,\dots,d}_{v\in B_n}$.

\reme In \cite{HWZ}, several possible variants of Rokhlin dimension have been exhibited for integer actions. Inserting $m=1$ in \ref{dimrok} yields the so-called Rokhlin dimension with single towers of \cite{HWZ}.

\rem \label{compact finite} Observe that the finite set $F\fin A$ in the definition above may be replaced by a compact set $K\subset A$. Suppose $K\subset A$ and $\eps>0$ is given. Then cover $K$ by balls of radius $\eps/3$. By compactness, we have
\[K\subset \bigcup_{a\in F} B_{\eps/3}(a)\quad\text{for a finite set}\; F\fin K.\]
Choose Rokhlin elements $(f_v^{(l)})_v^l$ for the pair $(F,\eps/3)$ and observe that these are in fact Rokhlin elements for the pair $(K,\eps)$.

\rem \label{decayfactors} For $n\in\IN$ and $j\in\set{-n+1,\dots,n}$, define
\[d_n(j) = 1-\frac{|j|}{n}.\]
Recall the notation $J_n= \set{-n+1,\dots,n}^m$.
For $m\in\IN$, denote $d^m_n: J_n \to [0,1]$ the map defined by
\[d^m_n(j_1,\dots,j_m) = \prod_{i=1}^m d_n(j_i).\]
For $a\in\set{0,1}^m$, let $s_a: J_n\to J_n$ be the bijection defined by
\[s_a[(j_i)_{i\leq m}] = (j_i+a_i\cdot n\mod 2n)_{i\leq m}.\]
One can easily see that for all $j\in\set{-n+1,\dots,n}$, one has 
\[1=d_n(j)+d_n(j+n\mod 2n).\]
Inductively, one can prove that then for all $v\in J_n$, one has
\[1=\sum_{a\in\set{0,1}^m} d^m_n(s_a(v)).\]

\rem \label{orthogonal} Let $A$ be a unital \cstar-algebra. Representing it on a Hilbert space $A\subset\CB(H)$, one can prove the following by application of the Cotlar-Stein Lemma (see \cite[chapter 7, section 2]{CoSt}). Let $J$ be a finite index set. For all $\delta>0$, there exists $\eta>0$ such that if $\set{b_j \;|\; j\in J}\subset A$ is a family of contractions with $\|b_i^*b_j\|\leq\eta$ for all $i\neq j$ in $J$, then
\[\|\sum_{j\in J} b_j\| \leq \delta+\max_{j\in J} \|b_j\| .\]
We now come to the main theorem of this section. Namely, we show that $\IZ^m$-actions of finite Rokhlin dimension preserve finite nuclear dimension of the underlying \cstar-algebra.

{\theorem[compare to {\cite[4.1]{HWZ}}] \label{Hauptsatz} Let $A$ be a unital \cstar-algebra and \linebreak $\alpha: \IZ^m\curvearrowright A$ be a group action via automorphisms. Then we have
\[\dimnuc(A\rtimes_{\alpha} \IZ^m)\leq 2^m(\dimnuc(A)+1)(\dimrokcyc(\alpha)+1)-1.\]}
\begin{proof}
We may assume that both $s=\dimnuc(A)$ and $d=\dimrokcyc(\alpha)$ are finite, or else the statement is trivial. 
Let $F\fin A\rtimes_{\alpha} \IZ^m$ a finite subset and $\delta>0$. Then we can find $N\in\IN$ and $F'\fin A\rtimes_{\alpha} \IZ^m$ such that $F=_\delta F'$ and
\[a=\sum_{v\in J_N} a(v)u_v\quad\text{for all}\; a\in F'.\]
Note that the $\set{ u_v | v\in\IZ^m}$ denote the canonical unitaries in $A\rtimes_\alpha\IZ^m$ implementing the action $\alpha$ on $A$. 
Let $\tilde{F}_1\fin A$ be the finite subset of all possible coefficients occuring in such a sum. Without loss of generality, we can assume that $\tilde{F}_1$ consists of contractions.
Let $\eps>0$ and $n\geq N$. As will be specified later, $n$ is a very large number compared to $N$ and $\eps$ very small compared to both $\delta$ and $n$. Let
\[\tilde{F} = \bigcup_{v\in J_n} \alpha^{-v}(\tilde{F}_1)\quad\fin\quad A.\]
\phantomsection\label{phipsi} Choose an $s$-decomposable c.p.~approximation $(\CF,\psi,\phi)$ for $\tilde{F}$ up to $\frac{\delta}{|J_n||J_N|}$, i.e. a finite dimensional \cstar-algebra $\CF=\CF^{(0)}\oplus\dots\oplus\CF^{(s)}$ and c.p.~maps $A\auf{\psi}\CF\auf{\phi} A$ such that $\psi$ is c.p.c., the maps $\phi^{(i)}=\phi|_{\CF^{(i)}}$ are c.p.c.~order zero and $\|x-(\phi\circ\psi)(x)\|\leq\frac{\delta}{|J_n||J_N|}$ for all $x\in\tilde{F}$.

Let $\psi_n = \id_{|J_n|}\otimes\psi: M_{|J_n|}\otimes A\to M_{|J_n|}\otimes\CF$ and $\phi_n = \id_{|J_n|}\otimes\phi: M_{|J_n|}\otimes\CF\to M_{|J_n|}\otimes A$ denote the amplifications of $\psi$ and $\phi$. Analogously write $\phi_n^{(i)} = \id_{|J_n|}\otimes\phi^{(i)}$.

Let $B_\CF$ be the closed unit ball in $\CF$ and define the compact set $K\subset A$ by
\[K = \bigcup_{v\in J_n} \bigcup_{i=0}^s \alpha^{-v}(\phi^{(i)}(B_\CF)). \]
Using \ref{compact finite}, we may choose positive contractions $(f^{(l)}_v)_{v\in B_{2n}}^{l=0,\dots,d}$ as in \ref{dimrok} for the triple $(K,\eps/2,2n)$. Moreover, we may arrange with a standard functional calculus argument that the relations (2)-(4) are also true for the square roots. By doing an index shift, we get elements $(f^{(l)}_v)_{v\in J_n}^{l=0,\dots,d}$ satisfying the analogous properties (1)-(4) (with the square roots in (2)-(4)) up to $\eps$ for the index set $J_n$ instead of $B_{2n}$.\vspace{3mm}

Let $A$ act faithfully on a Hilbert space $H$ and let $A\rtimes_{\alpha} \IZ^m$ be canonically embedded into $\CB(\ell^2(\IZ^m)\otimes H))$. Let $Q\in\CB(\ell^2(\IZ^m)\otimes H))$ be the projection onto the subspace $\ell^2(J_n)\otimes H$.
Then $x\mapsto QxQ$ defines a u.c.p. map \mbox{$\Psi: A\rtimes_{\alpha}\IZ^m \to M_{|J_n|}(A)$}. More specifically, we have for all $a\in A, v\in\IZ^m$ that
\[\begin{array}{cccc}
& \Psi(a u_v) &=& \displaystyle Q\Bigl[\sum_{w\in \IZ^m} e_{w,w-v}\otimes\alpha^{-w}(a)\Bigl]Q \\\\
= &  \displaystyle \sum_{w\in J_n: \atop w-v\in J_n} e_{w,w-v}\otimes\alpha^{-w}(a)
&=& \displaystyle \sum_{w\in J_n\cap (v+J_n)} e_{w,w-v}\otimes\alpha^{-w}(a).
\end{array}\]
Define the diagonal matrix $D\in M_{|J_n|}(\IC)$ by $D_{v,v} = d^m_n(v)$. Observe that if $v,w\in J_n$, then $|d^m_n(v)-d^m_n(w)|\leq \frac{m\|v-w\|_\infty}{n}$. If $n$ is large enough in comparison to $N$, we can ensure that
\[|\sqrt{d^m_n(v)}-\sqrt{d^m_n(w)}|\leq\delta/|J_N|\quad\text{for all}\; v,w\in J_n\;\text{with}\; v-w\in J_N.\]
It follows that for all $v\in J_N$ that
\[\begin{array}{lcl}
\multicolumn{3}{l}{\|[\sqrt{D},Qau_vQ]\|} \vspace{2mm}\\
\hspace{1cm}&=& \displaystyle \left\| \sum_{w\in J_n\cap (v+J_n)} \left( \sqrt{d^m_n(w)}-\sqrt{d^m_n(w-v)} \right) e_{w,w-v}\otimes \alpha^{-w}(a) \right\| \vspace{2mm}\\
&\leq& \max\set{ ~\left| \sqrt{d^m_n(w)}-\sqrt{d^m_n(w-v)} \right|~\bigl|~ w\in J_n\cap(v+J_n)}\cdot \|a\| \vspace{2mm}\\
&\leq& \frac{\delta}{|J_N|}\|a\|.
\end{array}\]
\phantomsection\label{mu}
Define the c.p.c.~map $\mu: A\rtimes_\alpha\IZ^m\to M_{|J_n|}(A)$ by $\mu(x)=\sqrt{D}\Psi(x)\sqrt{D}$. By the previous calculation, we have
\[\|\mu(au_v)-DQau_vQ\|\leq\delta/|J_N|\cdot\|a\|\quad\text{for all}\quad v\in J_N\;\text{and}\; a\in A.\]
Now let $p\in\set{0,1}^m$ and $l\in\set{0,\dots,d}$. Define maps $\sigma^{(l)}_p: M_{|J_n|}(A)\to A\rtimes_{\alpha} \IZ^m$ by $\sigma^{(l)}_p(e_{v,w}\otimes a) = f_{s_p(v)}^{(l) 1/2} u_{v} a u_{w}^* f_{s_p(w)}^{(l) 1/2}$ (see \ref{decayfactors} for notation). Note that these are indeed c.p.~since $\sigma^{(l)}_p(x)=v_{l,p} x v_{l,p}^*$ for the matrix $v_{l,p}\in M_{1, |J_n|}(A\rtimes_{\alpha} \IZ^m)$ defined by $v_{l,p}=(f^{(l) 1/2}_{s_p(w)} u_w)_{w\in J_n}$.

Now let $w_1,w_2,w_3,w_4\in J_n$ and $a,b\in K$. In the case $w_2=w_3$ we have
\[\begin{array}{lll}
\multicolumn{3}{l}{\sigma^{(l)}_0(e_{w_1,w_2}\otimes a)\sigma^{(l)}_0(e_{w_3,w_4}\otimes b)} \vspace{3mm}\\
\hspace{1cm}&=&  f_{w_1}^{(l) 1/2} u_{w_1} a u_{w_2}^* f_{w_2}^{(l) 1/2} \cdot f_{w_3}^{(l) 1/2} u_{w_3} b u_{w_4}^* f_{w_4}^{(l) 1/2} \vspace{3mm}\\
&=& f_{w_1}^{(l) 1/2} u_{w_1} a u_{w_2}^* f_{w_2}^{(l)} u_{w_2} b u_{w_4}^* f_{w_4}^{(l) 1/2} \vspace{3mm}\\
&\stackrel{(3)}{=}_{\eps}& f_{w_1}^{(l) 1/2} u_{w_1} a f_{0}^{(l)} b u_{w_4}^* f_{w_4)}^{(l) 1/2} \vspace{3mm}\\
&\stackrel{(4)}{=}_{\eps}& f_{w_1}^{(l) 1/2} u_{w_1} f_{0}^{(l)} a  b u_{w_4}^* f_{w_4}^{(l) 1/2} \vspace{3mm}\\
&\stackrel{(3)}{=}_{\eps}& f_{w_1}^{(l)}f_{w_1}^{(l) 1/2} u_{w_1} a  b u_{w_4}^* f_{w_4}^{(l) 1/2} \vspace{3mm}\\
&\stackrel{(2)}{=} & \hspace{-5mm}{}_{(|J_n|-1)\eps}\hspace{2mm}  f^{(l)} f_{w_1}^{(l) 1/2} u_{w_1} a  b u_{w_4}^* f_{w_4}^{(l) 1/2}.
\end{array}\]
For the last step of the above calculation, we denote $f^{(l)} = \sum_{w\in J_n} f^{(l)}_w$ for all $l=0,\dots, d$. In the case $w_2\neq w_3$ we have
\[ \begin{array}{lcl}
\multicolumn{3}{l}{\|\sigma^{(l)}_0(e_{w_1,w_2}\otimes a)\sigma^{(l)}_0(e_{w_3,w_4}\otimes b)\|} \\
\hspace{1cm}&=& \|f_{w_1}^{(l) 1/2} u_{w_1} a u_{w_2}^* f_{w_2}^{(l) 1/2} \cdot f_{w_3}^{(l) 1/2} u_{w_3} b u_{w_4}^* f_{w_4}^{(l) 1/2} \| \stackrel{(2)}{\leq} \eps.
\end{array}\]
So in any case, we get
\[\sigma^{(l)}_0(e_{w_1,w_2}\otimes a)\sigma^{(l)}_0(e_{w_3,w_4}\otimes b) =_{(|J_n|+2)\eps} f^{(l)}\sigma^{(l)}_0( (e_{w_1,w_2}\otimes a)(e_{w_3,w_4}\otimes b)).\]
It follows that
\[\sigma^{(l)}_0(a)\sigma^{(l)}_0(b) =_{(|J_n|+2)|J_n|^3\eps} f^{(l)}\sigma^{(l)}_0(ab) \]
for all $a,b\in M_{|J_n|}(K)$ and all $l=0,\dots,d$. The respective statements also hold for $\sigma^{(l)}_p$ in place of $\sigma^{(l)}_0$. So keeping in mind the definition of $K$, we have that $\sigma^{(l)}_p\circ\phi_n^{(i)}$ is order zero up to $(|J_n|+2)|J_n|^3\eps$ for all $i=0,\dots, s,\; p\in\set{0,1}^m$ and $l=0,\dots,d$.

\phantomsection \label{stern} Moreover, we have for all $l=0,\dots,d$,\; $v,w\in J_n$ and $a\in K$
\[\begin{array}{cll}
\sigma_0^{(l)}(e_{v,w}\otimes a) &=& f_v^{(l)1/2} u_v a u_{-w} f_w^{(l)1/2} \\\\
&\stackrel{(3)}{=}_{\eps}& f_v^{(l)1/2} u_v a f_0^{(l)1/2} u_{-w} \\\\
&\stackrel{(4)}{=}_{\eps}& f_v^{(l)1/2} u_v f_0^{(l)1/2} a u_{-w} \\\\
&\stackrel{(3)}{=}_{\eps}& f_v^{(l)} u_v a u_{-w}.
\end{array}\]
The respective statements also hold for $\sigma^{(l)}_p$ in place of $\sigma^{(l)}_0$. Let us denote
\[(\star)\qquad \sigma_p^{(l)}(e_{v,w}\otimes a) =_{3\eps} f_{s_p(v)}^{(l)} u_v a u_{-w}\quad\text{for all}\; l,\; v,w,\; p\;\text{and}\; a\in K.\]

Applying \ref{orthogonal} and (2), if $\eps$ is chosen small enough in comparison to $\delta$, we get for all $v\in J_N$ and $l=0,\dots,d$ that
\[\begin{array}{lll} \multicolumn{3}{l}{ \displaystyle \Bigl\|\sum_{w\in J_n\setminus (v+J_n)} d^m_n(w) \sum_{p\in\set{0,1}^m} f^{(l)}_{s_p(w)}\Bigl\| }
 \\\\
\hspace{3mm}&\leq & \displaystyle 2^m\cdot \max_p\; \Bigl\|\sum_{w\in J_n\setminus(v+J_n)} d^m_n(w) f^{(l)}_{s_p(w)}\Bigl\| \\\\
&\stackrel{\ref{orthogonal}}{\leq}& 2^m\cdot\left[\delta/|J_N|+\max\set{ \|d^m_n(w)f_{s_p(w)}^{(l)}\| \;|\; w\in J_n\setminus (v+J_n),\; p\in\set{0,1}^m}\right] \\\\
&\leq& 2^m\cdot\left[\delta/|J_N|+\max\set{ d^m_n(w) \;|\; w\in J_n\setminus (v+J_n)}\right] \\\\
&\leq& 2^m\cdot\left[\delta/|J_N|+(1-\frac{n-N}{n})^m\right].
\end{array}\]  
We may assume that $n$ is large enough in comparison to $N$ such that the right side is less or equal $2^{m+1}\delta/|J_N|$. Observe that we have for all $l$
\[\sum_{w\in J_n} f_w^{(l)} = \sum_{w\in J_n}  d^m_n(w) \sum_{p\in\set{0,1}^m} f^{(l)}_{s_p(w)}.\]
\phantomsection\label{summeeins} It follows for all $v\in J_N$ that
\[\begin{array}{lll} \multicolumn{3}{l}{\displaystyle \Bigl\|\eins_A-\sum_{l=0}^d\sum_{w\in J_n\cap (v+J_n)} d^m_n(w) \sum_{p\in\set{0,1}^m} f^{(l)}_{s_p(w)} \Bigl\| } \\\\
\hspace{1cm}&\stackrel{(1)}{\leq}& \displaystyle \eps+ \Bigl\|\sum_{l=0}^d \sum_{w\in J_n\setminus (v+J_n)} d^m_n(w) \sum_{p\in\set{0,1}^m} f^{(l)}_{s_p(w)} \Bigl\| \\\\
&\leq & \eps+(d+1)2^{m+1}\delta/|J_N| \quad\leq\quad \frac{2^{m+2}(d+1)\delta}{|J_N|}.
\end{array}\]

Now let $l\in\set{0,\dots,d}, a\in\tilde{F}$ and $v\in J_N$. Denote $\sigma:= \sum_{l=0}^d\sum_{p\in\set{0,1}^m} \sigma^{(l)}_p$. We have
\begin{longtable}{ll}
\multicolumn{2}{l}{$\sigma\circ\phi_n\circ\psi_n\circ\mu(au_v)$} \\\\
$=$ & \hspace{-25mm}$\displaystyle \sum_{l=0}^d\sum_{p\in\set{0,1}^m}\sum_{i=0}^s (\sigma^{(l)}_p\circ\phi^{(i)}_n\circ\psi_n\circ\mu)(au_v)$ \\\\
$\hyperref[mu]{=}_{2^m(d+1)(s+1)\delta/|J_N|}$ & $\displaystyle \sum_{l=0}^d\sum_p\sum_{i=0}^s (\sigma^{(l)}_p\circ\phi_n^{(i)}\circ\psi_n)(DQau_vQ)$ \\\\
$=$ & \hspace{-25mm}$\displaystyle \sum_{l=0}^d\sum_p \sum_{i=0}^s (\sigma_p^{(l)}\circ\phi_n^{(i)}\circ\psi_n) \left( \sum_{w\in J_n\cap (v+J_n)} d^m_n(w)\cdot e_{w,w-v}\otimes\alpha^{-w}(a) \right)$ \\\\
$=$ & \hspace{-25mm}$\displaystyle \sum_{l=0}^d\sum_p \sum_{i=0}^s \sigma_p^{(l)} \left( \sum_{w\in J_n\cap (v+J_n)} d^m_n(w)\cdot e_{w,w-v}\otimes (\phi^{(i)}\circ\psi)(\alpha^{-w}(a)) \right)$ \\\\
$\hyperref[stern]{\stackrel{(\star)}{=}}_{2^m|J_n|(d+1)(s+1)\cdot 3\eps}$ 
& \hspace{-2mm}$\displaystyle \sum_{l=0}^d\sum_p\sum_{i=0}^s \sum_{w} d^m_n(w)f_{s_p(w)}^{(l)}\cdot u_w [(\phi^{(i)}\circ\psi)(\alpha^{-w}(a))] u_{v-w}$ \\\\
$=$ & \hspace{-25mm}$\displaystyle \sum_{l=0}^d\sum_p \sum_{w} d^m_n(w)f_{s_p(w)}^{(l)}\cdot u_w [(\phi\circ\psi)(\alpha^{-w}(a))] u_{v-w}$ \\\\
$\hyperref[phipsi]{=}_{2^m(d+1)|J_n|\cdot \frac{\delta}{|J_n||J_N|}}$ & \hspace{-1mm}$\displaystyle \sum_{l=0}^d\sum_p \sum_{w} d^m_n(w)f_{s_p(w)}^{(l)}\cdot u_w \alpha^{-w}(a) u_{v-w}$ \\\\
$=$ & \hspace{-25mm}$\displaystyle \left(\sum_{l=0}^d\sum_p \sum_{w} d^m_n(w)f_{s_p(w)}^{(l)}\right)\cdot a u_v$ \\\\
$\hyperref[summeeins]{=}_{2^{m+2}(d+1)\delta/|J_N|}$ & $au_v$.
\end{longtable}
Assuming that $3|J_n|\eps\leq\delta/|J_N|$, it follows for all $a\in \tilde{F}$ and $v\in J_N$ that
\[\|au_v-(\sigma\circ\phi_n\circ\psi_n\circ\mu)(au_v)\|\leq 2^{m+3}(d+1)(s+1)\cdot \frac{\delta}{|J_N|}.\]
Hence for $x\in F'$ we get
\[\|x-(\sigma\circ\phi_n\circ\psi_n\circ\mu)(x)\|\leq 2^{m+3}(d+1)(s+1)\delta.\] 
Observe that $\|\sigma\|\leq 2^m(d+2)$ for small enough $\eps$. Hence $\|\sigma\circ\phi_n\|\leq 2^m(d+2)(s+1)$. Thus, we get for all $x\in F$
\[\begin{array}{lcl}
\multicolumn{3}{l}{\|x-(\sigma\circ\phi_n\circ\psi_n\circ\mu)(x)\|} \vspace{2mm}\\ 
\hspace{1cm}&\stackrel{F\approx_\delta F' \atop \exists y\in F}{\leq}& \delta+2^m(d+2)(s+1)\delta+\|y-(\sigma\circ\phi_n\circ\psi_n\circ\mu)(y)\| \vspace{2mm}\\
&\leq & 2^{m+4}(d+1)(s+1)\delta.
\end{array}\]
Now let us recall what we got. We have contructed a c.p.~approximation 
\[(M_{|J_n|}(\CF)^{2^m(d+1)},\; \psi_n\circ\mu,\; \sigma\circ\phi_n)\]
 of tolerance $2^{m+4}(d+1)(s+1)\delta$ on $F$, where the map
\[\sigma\circ\phi_n = \sum_{i=0}^s\sum_{l=0}^d\sum_{p\in\set{0,1}^m} \sigma^{(l)}_p\circ\phi_n^{(i)}\]
is a direct sum of $2^m(d+1)(s+1)$ c.p.c.~maps of order zero up to $(|J_n|+2)|J_n|^3\eps$. Since we can always choose $\eps$ as small as we want in relation to $n$, the maps $\sigma_p^{(l)}\circ\phi_n^{(i)}$ can be chosen to be order zero up to $\eta$ for any given $\eta>0$.\vspace{3mm}

Since $m,d,s$ are constants and $F\fin A\rtimes_\alpha\IZ^m$ and $\delta>0$ were arbitrary, it follows from \ref{orderdelta} that
\[\dimnuc(A\rtimes_\alpha\IZ^m)\leq 2^m(d+1)(s+1)-1,\]
which is what we wanted to show.
\end{proof}


\section{The Rokhlin dimension for topological $\IZ^m$-actions}
\noindent
Since we are (at least here) mainly interested in \cstar-dynamical systems arising from topological dynamical systems, it is natural to ask whether there exists a variant of Rokhlin dimension for these actions that is phrased in topological terms. The aim of this section is to give such a variant and to show that finite Rokhlin dimension in the topological sense implies finite Rokhlin dimension in the \cstar-algebraic sense.

The definition of the topological variant is due to Wilhelm Winter.
\defi \label{dimrok top} Let $(X,\alpha,\IZ^m)$ be a topological dynamical system. We say that $\alpha$ has Rokhlin dimension $d$, and write $\dimrok(\alpha)=d$, if $d$ is the smallest natural number with the following property:

For all $n\in\IN$, there exists a family of open sets 
\[\CR=\set{ U^{(l)}_v \;|\; l=0,\dots,d,\; v\in B_n}\]
in $X$ such that
\begin{itemize}
\item $U^{(l)}_v=\alpha^v(U^{(l)}_0)$ for all $l=0,\dots,d$ and $v\in B_n$.
\item For all $l$, the sets $\set{ \quer{U}^{(l)}_v \;|\; v\in B_n}$ are pairwise disjoint.
\item $\CR$ is an open cover of $X$.
\end{itemize} 
If there is no such $d$, we write $\dimrok(\alpha)=\infty$. In this context, we call $\CR$ a Rokhlin cover.

{\lemma[cf. {\cite[Proposition 2.7]{HWZ}}] \label{schwacher dimrok} Let $\alpha: \IZ^m\curvearrowright X$ be a free continuous action on a compact metric space. Let $\bar{\alpha}: \IZ^m\curvearrowright\CC(X)$ denote the induced \cstar-algebraic action by $\alpha$. Let $d\in\IN$ be a natural number with the properties of \ref{dimrok}, with the exception that property (1) is replaced by the condition
\[(1') \quad \sum_{l=0}^d \sum_{v\in B_n} f^{(l)}_v \geq \eins_A.\]
Then $\dimrok(\bar{\alpha})\leq d$.}
\begin{proof}
Let $\eps,\delta>0$ and $n\in\IN$ be given. Choose positive contractions 
\[\set{ h^{(l)}_v \;|\; l=0,\dots,d,\; v\in B_n}\]
satisfying properties (1') and (2)-(4) for $\delta$. If $\delta$ is small enough, we can ensure that
\[1\leq\sum_{l=0}^d \sum_{v\in B_n} h_v^{(l)} \leq d+2.\]
Denote the function in the middle by $S$. We have 
\[\|S-S\circ\alpha^{-v}\|\leq |B_n|(d+1)\delta\quad\text{for all}\; v\in B_n.\]
By approximating the function $[1,d+2]\to [0,1],\; t\mapsto t^{-1}$ by polynomials, one gets that
\[\|S^{-1}-S^{-1}\circ\alpha^{-v}\|\leq\eps/2 \quad\text{for all}\; v\in B_n,\]
if $\delta$ is small enough in relation to $\eps$ and $n$. Now define $f^{(l)}_v=S^{-1}\cdot h^{(l)}_v$ for all $l=0,\dots, d$ and $v\in B_n$.
Observe that these elements satisfy the following properties:
\begin{itemize}
\item The $f_v^{(l)}$ are positive contractions.
\item $\sum_{l=0}^d \sum_{v\in B_n} f^{(l)}_v = \eins_A$.
\item Property (4) holds trivially.
\item Property (2) still holds for $\delta$.
\item As for property (3), we get for all $v,w\in B_n$
\[\begin{array}{ccl}
\|\bar{\alpha}^v(f^{(l)}_w)-f^{(l)}_{v+w}\| &=& \| (S^{-1}\circ\alpha^{-v})\cdot\bar{\alpha}^v(h^{(l)}_w)-S^{-1}\cdot h^{(l)}_{v+w}\| \\\\
&\leq & \eps/2 + \|S^{-1}\|\cdot\|\bar{\alpha}^v(h^{(l)}_w)-h^{(l)}_{v+w}\| \\\\
&\leq & \eps/2+\delta.
\end{array}\]
\end{itemize}
We see that if $\delta$ is chosen small enough in relation to $\eps$ and $n$, then the new functions $f^{(l)}_v$ satisfy the relations (1)-(4) for $\eps$.
\end{proof}

\rem A slightly more general statement is true by almost the same proof. Namely, if $A$ is a unital \cstar-algebra with an action $\alpha: \IZ^m\curvearrowright A$, let $d$ be a natural number with the properties (2)-(5) of \ref{dimrok} and the property (1'). Then $\dimrokcycc(\alpha)\leq d$. In other words, in the above theorem, one only needs commuting towers instead of a commutative \cstar-algebra.

{\prop \label{cyc und notcyc} Let $\alpha: \IZ^m\curvearrowright X$ be a continuous group action on a compact metric space. Let $\bar{\alpha}: \IZ^m\curvearrowright\CC(X)$ denote the induced \cstar-algebraic action by $\alpha$. Then
\[\dimrokcyc(\bar{\alpha}) \leq 2^m(\dimrok(\alpha)+1)-1.\]}
\begin{proof}
We may certainly assume that the right side is finite. Let $d=\dimrok(\alpha)$. For $N\in\IN$, let $\set{ U^{(l)}_{v}\;|\; l=0,\dots,d,\; v\in B_{2N} }$ be a Rokhlin cover. Note that by doing an index shift, we may as well assume the form $\set{ U_{v}^{(l)} \;|\; l=0,\dots,d,\; v\in J_N}$ with the corresponding relations, which is notationally more convenient in this proof.

Fix some $l$ and find $\delta>0$ small enough such that the sets $\set{ V_v^{(l)} \;|\; v\in J_N }$ are still pairwise disjoint for $V_v^{(l)} = \alpha^v(B_\delta(U_0^{(l)}))$. Let $h\in\CC_0(V_0^{(l)})$ be a function with $h|_{U_0^{(l)}}=1$. For all $v\in J_N$, define $g^{(l)}_v\in\CC_0(V_v^{(l)})$ via $g^{(l)}_v(x) = h(\alpha^{-v}(x)).$

Observe the following properties:
\begin{itemize}
\item $g^{(l)}_v = g^{(l)}_0\circ\alpha^{-v}$ for all $v\in J_N$.
\item $g^{(l)}_v \cdot g^{(l)}_w = 0$ for $v\neq w$ in $J_N$.
\item $g^{(l)}_v$ is constantly 1 on $U^{(l)}_v$ for all $v\in J_N$.
\end{itemize}

Now let $L\in\IN$ and $\eps>0$. Choose $n$ large enough such that $\frac{2}{n}\leq\eps$. Choose a Rokhlin cover $\set{U_v^{(l)} \;|\; l=0,\dots,d,\; v\in J_{N}}$ for $N=4Ln$ and choose the functions $\set{ g^{(l)}_v \;|\; l=0,\dots,d,\; v\in J_{N}}$ as above. For all $l$, define functions $(f^{(l)}_v)_{v\in B_{L}}$ via
\[f^{(l)}_v(x) = \begin{cases} g^{(l)}_w(x) &,\quad \text{if}\; x\in V_w^{(l)}\;\text{for}\; \|w\|_\infty\leq 2Ln\\
& \hspace{5.5mm}\text{and}\; w=v\mod L\IZ^m \\\\
\frac{3Ln-\|w\|_\infty}{Ln}\cdot g^{(l)}_w(x) &,\quad \text{if}\; x\in V_w^{(l)}\;\text{for}\; 2Ln<\|w\|_\infty\leq 3Ln\\
& \hspace{5.5mm} \text{and}\; w=v\mod L\IZ^m \\\\
0 &,\quad\text{elsewhere}.
 \end{cases}\]
Now the properties of $g^{(l)}_v$ ensure that
\begin{itemize}
\item $f^{(l)}_v\cdot f^{(l)}_w=0$ for $v\neq w$ in $B_L$. 
\item $\sum_{v\in B_L} f^{(l)}_v$ is constantly 1 on $\bigcup_{w\in I_{2Ln}} U_w^{(l)}$.
\item For $\|w\|_\infty=L$, we have $\|f_0^{(l)}-f_0^{(l)}\circ\alpha^{-w}\|\leq \frac{1}{n}$.
\item For $v\in B_L$, we have $\|f_v^{(l)}-f_0^{(l)}\circ\alpha^{-v}\|\leq \frac{1}{n}$.
\item Hence $\|f_w^{(l)}\circ\alpha^{-v}-f^{(l)}_{(v+w)\mod L\IZ^m}\|\leq\frac{2}{n}\leq\eps$ for all $v,w\in B_L$.
\end{itemize}
Now choose $\set{ a_j \;|\; j=1,\dots,2^m}\subset\IZ^m$ such that 
\[J_{N}=J_{4Ln}=\bigcup_{j=1}^{2^m} a_j+J_{2Ln}.\]
For $l=0,\dots,d,\;j=1,\dots,2^m$ and $v\in B_L$, define $f^{(l,j)}_v = f^{(l)}_v\circ\alpha^{-a_j}$. We have established that for all $l$ and $j$, the functions $(f^{(l,j)}_v)_{v\in B_L}$ satisfy the relations (2)-(4) of \ref{dimrok}. Furthermore, $\sum_{v\in B_L} f^{(l,j)}_v$ is constantly 1 on $\bigcup_{w\in (a_j+J_{2Ln})} U_{w}^{(l)}$, so the choice of the $a_j$ ensures that $\sum_{j=1}^{2^m}\sum_{v\in B_L} f_v^{(l,j)}\geq 1$ on $\bigcup_{w\in J_{N}} U^{(l)}_w$, hence 
\[\sum_{l=0}^d \sum_{j=1}^{2^m}\sum_{v\in B_L} f_v^{(l,j)}\geq 1.\]
We see that the family $\set{ f^{(l,j)}_v\;|\; l,j,v}$ has $2^m\cdot(d+1)$ upper indices, so we are done by Lemma \ref{schwacher dimrok}.
\end{proof}
We see now that it suffices to study the topological variant of Rokhlin dimension of a topological $\IZ^m$-action in order to study the \cstar-algebraic counterpart. A result of \cite{Gutman} enables an easy proof of finite Rokhlin dimension for $m=1$ in this vein:

\begin{samepage}\theorem[see {\cite[theorem 6.1]{Gutman}} with proof.]\label{Gutman} Let $X$ be a compact metric space of finite covering dimension $d$ and $\phi: X\to X$ an aperiodic homeomorphism. For all $n\in\IN$, we can find an open set $Z\subset X$ such that
\begin{itemize}
\item $\quer{Z}\cap\phi^j(\quer{Z})=\emptyset$ for all $j=1,\dots,n-1$.
\item $X=\bigcup_{j=0}^{2(d+1)n-1} \phi^j(Z).$
\end{itemize}\end{samepage}

Such an open set is called an $n$-marker. We will generalize this notion in Section 4 for arbitrary group actions.

{\cor \label{dimrok endlich Z} Let $X$ be a compact metric space and $\phi: X\to X$ an aperiodic homeomorphism. Then we have 
\[\dimrok(\phi)\leq 2(\dim(X)+1)-1.\]}
\begin{proof}
We may assume that $X$ has finite covering dimension $d$, or else there is nothing to show. Let $n\in\IN$ be given. Find an open set $Z\subset X$ like in \ref{Gutman} for this number. For $j=0,\dots,n-1$ and $l=0,\dots,2d+1$ set
\[U_j^{(l)}=\phi^{ln+j}(Z).\]
Then it is obvious that $\CR=\set{ U_j^{(l)} \;|\; l=0,\dots,2d+1,\; j=0,\dots,n-1}$ is a Rokhlin cover. That is, we have
\begin{itemize}
\item $U^{(l)}_j=\phi^j(U^{(l)}_0)$ for all $l$ and $j$.
\item For all $l$, the sets $\set{ \quer{U}^{(l)}_j \;|\; j=0,\dots,n-1}$ are pairwise disjoint.
\item $\CR$ is an open cover of $X$.
\end{itemize} 
Hence $\dimrok(\phi)\leq 2d+1$.
\end{proof}


\section{The topological small boundary property}
\noindent
In this section, we define a technical condition that we name the (bounded) topological small boundary property. Weaker versions of this were considered by Lindenstrauss in \cite{Lind}, by Gutman in \cite{Gutman2} and \cite{Gutman3} and had connections to a dynamical system having mean dimension zero. To the author's knowledge, Phillips was the first one to link the topological small boundary property to strict comparison for the crossed product, thereby linking it to \cstar-classification as well (see \cite{Phi}). Note also \cite{Buck}, where the topological small boundary property is shown to imply a purely dynamical analogue of strict comparison. 

It will turn out that we can in fact assume a stronger bounded variant, whenever we have a finite dimensional underlying space.

\defi \label{disjointness} Let $(X,\alpha,G)$ be a topological dynamical system. Let $M\subset G$ be a subset and $k\in\IN$ be some natural number. We say that a set $E\subset X$ is $(M,k)$-disjoint, if for all distinct elements $\gamma(0),\dots,\gamma(k)\in M$ we have
\[\alpha^{\gamma(0)}(E)\cap\dots\cap\alpha^{\gamma(k)}(E)=\emptyset.\]
We call $E$ toplogically $\alpha$-small if $E$ is $(G,k)$-disjoint for some $k$. We call the smallest such $k$ the (topological) smallness constant of $E$.

\defi[cf. \cite{Lind}] \label{SBP}  Let $G$ be a countably infinite group. A topological dynamical system $(X,\alpha,G)$ has the topological small boundary property, if whenever $A\subset X$ is closed and $U\supset A$ is open, we can find $A\subset V\subset U$ open such that $\del V$ is topologically $\alpha$-small.

\phantomsection\label{TSBPd} If we can arrange that each such $\del V$ has a smallness constant bounded uniformly by a number $d$, we say that $(X,\alpha,G)$ has the bounded topological small boundary property with respect to $d$, abbreviated (TSBP $\leq d$).\vspace{3mm}\\
The main goal of this section is to prove that free actions on finite dimensional spaces have this property. It is important to note that the case $G=\IZ$ has been treated by Lindenstrauss in \cite{Lind}. He has shown that, if $\phi: X\to X$ is an aperiodic homeomorphism, then $(X,\phi)$ has the bounded topological small boundary property with respect to $\dim(X)$. We merely give a modification of his proof to obtain essentially the same result for free countable group actions.\vspace{3mm}\\
We would like to quote some well-known facts about properties of covering dimension, since using these will be crucial in some key steps. These statements come up in \cite{Lind}, but for a detailed treatment the reader is referred to \cite{Engel}.
All spaces in question are assumed to be separable metric spaces.
\begin{itemize}
\item[D1]\phantomsection\label{D1} $A\subset B$ implies $\dim(A)\leq \dim(B)$.
\item[D2]\phantomsection\label{D2} If $\set{B_i}_{i\in\IN}$ is a family of closed sets in $A$ with $\dim(B_i)\leq k$, then $\dim(\bigcup B_i)\leq k$.
\item[D3]\phantomsection\label{D3} Let $E \subset A$ be a zero dimensional subspace and $x\in U\subset A$ a point with an open neighbourhood. There exists some open $U'\subset A$ with $x\in U'\subset U$ such that $\del U' \cap E=\emptyset$.
\item[D4]\phantomsection\label{D4} If $A\neq\emptyset$, there exists $E\subset A$, which is the countable union of zero dimensional closed subsets in $A$, such that $\dim(A\setminus E)=\dim(A)-1$.
\item[D5]\phantomsection\label{D5} Any countable union of $k$-dimensional $F_\sigma$-sets is a $k$-dimensional $F_\sigma$-set.
\end{itemize}

{\lemma \label{Startlemma} Let $X$ be a compact metric space and $U, V$ open sets with $\quer{U}\subset V$. Let $E\subset X$ be a zero dimensional subspace. There exists an open set $U'$ with $\quer{U}\subset U'\subset V$ such that $\del U'\cap E=\emptyset.$}
\begin{proof}
Clearly $\del U$ is compact. For $x\in\del U$, apply (D3) and find open neighbourhoods $x\in B_x\subset V$ such that $\del B_x \cap E=\emptyset$. Choose a finite cover $\del U\subset \bigcup_{i=1}^M B_i$ of such neighbourhoods and set $U'=U\cup\bigcup_{i=1}^M B_i$. It is now easy to see that $\del U' \subset \bigcup_{i=1}^M \del B_i$, so we have indeed $\del U' \cap E=\emptyset$.
\end{proof}

\defi Let $X$ be a compact metric space, $G$ a countable group and $\alpha: G\curvearrowright X$ a continuous group action. For any $g\in G\setminus\set{e}$, set
\[X^g = \set{x\in X\;|\; \alpha^g(x)=x}.\]
Moreover, let
\[X^G = \bigcup_{g\in G\setminus\set{e}} X^g.\]
We call $X^g$ the $g$-fixed point set. If $\dim(X^g)\leq 0$ for all $g\neq e$, we say that $\alpha$ has at most zero dimensional fixed point sets. By \hyperref[D5]{D5}, this is equivalent to $\dim(X^G)\leq 0$. Obviously, free actions fall under this category.

\defi \label{genposition} Let $X$ be a compact metric space of finite covering dimension $n$. A family
$\CB$ of subsets in $X$ is in general position, if for all finite subsets $S\fin\CB$ we have
\[\dim(\bigcap S) \leq \max(-1,n-|S|).\]

{\lemma \label{general position} Let $X$ be a compact metric space of finite covering dimension $n$ together with a continuous group action $\alpha: G\curvearrowright X$ that has at most zero dimensional fixed point sets.

Let $U,V\subset X$ be open sets with $\quer{U}\subset V$. For any finite subset $F\fin G$, there exists an open set $U'$ with $U\subset U'\subset \quer{U}'\subset V$ such that the family $\set{\alpha^\gamma(\del U')}_{\gamma\in F}$ is in general position and $\del U'\cap X^G=\emptyset$.}
\begin{proof}
We prove this by induction in the variable $k=|F|$. The assertion holds for $k=1$ by application of \ref{Startlemma}. Now assume that the assertion holds for some natural number $k$. We show that it also holds for $k+1$.

Let $F=\set{\gamma(0),\dots,\gamma(k)}$ be a set of cardinality $k+1$ in $G$. Using the induction hypothethis, there exists an open $A_0$ with $U\subset A_0\subset \quer{A_0}\subset V$, such that the collection $\set{\alpha^{\gamma(0)}(\del A_0),\dots,\alpha^{\gamma(k-1)}(\del A_0)}$ is in general position and $\del A_0\cap X^G=\emptyset$.

Because of $\quer{A_0}\subset V$ and the fact that $\del A_0\subset X\setminus X^G$, for all $x\in\del A_0$ we can find $\eta(x)>0$ such that $\quer{B}_{\eta(x)}(x)\subset V$ and such that $(\alpha^{\gamma(j)}(B_{\eta(x)}(x)))_{j=0,\dots,k}$ are pairwise disjoint. Denote $\hat{B_x} = B_{\eta(x)}(x)$ and $B_x=B_{\eta(x)/2}(x)$.

Find some finite subcover $\del A_0 \subset \bigcup_{i=1}^M B_i$. We will now construct open sets $A_i$ for $i=0,\dots,M$ ($A_0$ is already defined) with the following properties:
\begin{itemize}
\item[(0)] $\del A_i \cap X^G=\emptyset$.
\item[(1)] $\quer{A_i} \subset A_0\cup\bigcup_{j=1}^M B_j$.
\item[(2)] $A_{i}\subset A_{i+1}\subset A_{i}\cup \hat{B}_{i+1}$.
\item[(3)] The collection 
\[\CA_i = \set{ \alpha^{\gamma(j)}(\del A_i)}_{j<k}\cup
\Bigl\{\alpha^{\gamma(k)}(\del A_i \cap \bigcup_{j=1}^i B_j) \Bigl\}\]
 is in general position. 
\end{itemize}
Since by (1) we clearly have $\quer{A}_M\subset V$, combining this with $(3)$ implies that the set $U'=A_M$ has the desired property once we have done this construction. It remains to show how to construct the sets $A_i$. 

So suppose that the set $A_i$ has already been defined for $i<M$. According to \hyperref[D4]{D4}, for all nonempty subsets $S\subset\CA_i$, there exists a zero dimensional $F_\sigma$-set $E_S\subset\bigcap S$ such that $\dim(\bigcap S\setminus E_S) = \dim(\bigcap S)-1.$
Define
\[E:= \bigcup_{\emptyset\neq S\subset \CA_i \atop 0\leq j\leq k} \alpha^{-\gamma(j)}(E_S).\]
By \hyperref[D5]{D5}, $E$ is a zero dimensional $F_\sigma$-set. Use \ref{Startlemma} to find an open set $W$ such that
\[(\star)\qquad \quer{A_i\cap B_{i+1}} \subset W\subset\quer{W}\subset \hat{B}_{i+1}\cap (A_0\cup\bigcup_{j=1}^M B_j)\quad\text{and}\quad \del W \cap (E\cup X^G)=\emptyset.\]
Now set $A_{i+1}:= A_i\cup W$. This clearly satisfies the properties (0),(1) and (2). To show (3), let $\emptyset\neq S=\set{S_1,\dots,S_m}\subset \CA_{i+1}$ correspond to some subset $\set{j_1,j_2,\dots,j_m}\subset\set{0,\dots,k}$. Note that since $\del A_{i+1}\subset\del A_i\cup\del W$, we have either
\[S_l=\alpha^{\gamma(j_l)}(\del A_{i+1}) \subset \alpha^{\gamma(j_l)}(\del A_i)\cup 
\alpha^{\gamma(j_l)}(\del W) =: S_l^0\cup S_l^1\qquad(\text{if}\; j_l\neq k)\]
or
\[\begin{array}{ccc} 
S_l &=&\displaystyle \alpha^{\gamma(j_l)}(\del A_{i+1}\cap\bigcup_{j=1}^{i+1} B_i)\\\\
&\subset& \displaystyle \alpha^{\gamma(j_l)}((\del A_i\setminus W)\cap\bigcup_{j=1}^{i+1} B_i)\cup \alpha^{\gamma(j_l)}(\del W) \\\\
&\stackrel{(\star)}{\subset}& \displaystyle \alpha^{\gamma(j_l)}(\del A_i\cap\bigcup_{j=1}^i B_j)\cup\alpha^{\gamma(j_l)}(\del W) \\\\
 &=:& S_l^0\cup S_l^1\qquad\qquad(\text{if}\; j_l=k).
\end{array}\]
It follows that
\[\bigcap S \subset \bigcup_{a\in\set{0,1}^m} \left( \bigcap_{l=1}^m S_l^{a_l} \right).\]
Since $\quer{W}\subset \hat{B}_{i+1}$, our choice of $\hat{B}_{i+1}$ implies that the sets $S_l^1$ are pairwise disjoint. So it suffices to consider the case $a=(0,\dots,0)$ and, since we can change the order without loss of generality, the case $a=(1,0,\dots,0)$. For $a=(0,\dots,0)$, note that $\set{S_1^0,\dots,S_m^0}$ is a subset of $\CA_i$, so we already have
\[\dim\Bigl(\bigcap_{l=1}^m S_l^0\Bigl) \quad\leq\quad \max(-1,n-m).\]
For $a=(1,0,\dots,0)$, define $\hat{S}=\set{S_2^0,\dots,S_m^0}$. This is a subset of $\CA_i$, hence we know that it is in general position. Moreover, considering our choice of the set $E_{\hat{S}}$, recall that
\[\dim(\bigcap \hat{S}\setminus E_{\hat{S}}) \leq\dim(\bigcap\hat{S})-1 \leq \max(-1,n-(m-1))-1 \leq \max(-1,n-m).\]
By the choice of $W$ we know that $\del W\cap E=\emptyset$. Since $\alpha^{-\gamma(j_1)}(E_{\hat{S}})\subset E$, this implies $E_{\hat{S}}\cap\alpha^{\gamma(j_1)}(\del W)=\emptyset$. In particular, it follows that
\[S_1^1\cap\bigcap_{l=2}^m S_l^0 = \alpha^{\gamma(j_1)}(\del W)\cap\bigcap \hat{S} \subset \bigcap\hat{S}\setminus E_{\hat{S}}.\]
Therefore we have established
\[\dim(S_1^1\cap\bigcap_{l=2}^m S_l^0)\leq\max(-1,n-m).\]
If we combine these inequalities with \hyperref[D2]{D2}, it follows that we have $\dim(\bigcap S) \leq \max(-1,n-m)$ as well. So $\CA_{i+1}$ is in general position and we are done.
\end{proof} 

{\lemma \label{disjointness neighbourhood} Let $(X,\alpha,G)$ as before, $F\fin G$ a finite subset and $n\in\IN$ a natural number. If a closed subset $E\subset X$ is $(F,n)$-disjoint, there exists an open neighbourhood $V$ of $E$ such that $\quer{V}$ is $(F,n)$-disjoint.}
\begin{proof}
Note that for all $S\subset F$ with $n=|S|$, we have
\[\emptyset=\bigcap_{\gamma\in S} \alpha^\gamma(E)= \bigcap_{\gamma\in S} \alpha^\gamma\Bigl( \bigcap_{\eps>0} \quer{B}_\eps(E) \Bigl)
= \bigcap_{\eps>0} \left( \bigcap_{\gamma\in S} \alpha^\gamma(\quer{B}_\eps(E)) \right)\]
By compactness, there must exist some $\eps(S)>0$ such that 
\[\displaystyle \bigcap_{\gamma\in S} \alpha^\gamma(\quer{B}_{\eps(S)}(E))=\emptyset.\]
If we set $\eps=\min\set{ \eps(S) | S\subset F, n=|S|}$, the open set $V=B_\eps(E)$ clearly does the trick.
\end{proof}

{\theorem \label{TSBP} Let $X$ be a compact metric space with finite covering dimension $d$. Let $G$ be a countably infinite group with a continuous action $\alpha: G\curvearrowright X$ that has at most zero dimensional fixed point sets.
Then $(X,\alpha,G)$ satisfies \hyperref[TSBPd]{\em(TSBP $\leq d$)}. More precisely:

Let $U,V\subset X$ be open sets with $\quer{U}\subset V$. Then there exists an open set $U'$ with $U\subset U'\subset V$ such that $\del U'$ is $(G,d)$-disjoint.}
\begin{proof}
Since $G$ is countable, choose an increasing sequence $F_1\subset F_2\subset F_3\subset\dots\fin G$ of finite sets such that $G=\bigcup_{k=1}^\infty F_k$.
We will construct a sequence of open sets $\set{U_k}_{k\in\IN}, \set{V_k}_{k\in\IN}$ with the following properties for all $k$:
\begin{itemize}
\item[(1)] $U\subset U_k\subset U_{k+1}\subset V$.
\item[(2)] $\quer{V}_{k+1}\subset V_k$.
\item[(3)] $\quer{U}_{k+1}\subset U_k\cup V_k$.
\item[(4)] $\quer{V}_k$ is $(F_k,d)$-disjoint.
\item[(5)] $\del U_k \subset V_k$.
\end{itemize}
Set $U_0=U, V_0=V$. Apply \ref{general position} to find $U_1$ such that $U\subset U_1\subset \quer{U}_1\subset V$ and $\del U_1$ is $(F_1,d)$-disjoint. Apply \ref{disjointness neighbourhood} to find an open neighbourhood $V_1$ of $\del U_1$ such that $\quer{V}_1\subset V$ and $\quer{V}_1$ is $(F_1,d)$-disjoint. Clearly these sets satisfy (1)-(5) thus far.
Suppose that the sets $U_k, V_k$ have been defined for some $k$. Apply \ref{general position} to find an open set $U_{k+1}$ such that $U_k\subset U_{k+1}\subset \quer{U}_{k+1}\subset U_k\cup V_k$ and $\del U_{k+1}$ is $(F_{k+1},d)$-disjoint. Since $V_k$ is an open neighbourhood of $\del U_{k+1}$, we can find an open neighbourhood $V_{k+1}$ of $\del U_{k+1}$ such that $\quer{V}_{k+1}\subset V_k$ is $(F_{k+1},d)$-disjoint. It is easy to see that these new sets satisfy properties (1)-(5) again.

Now set $U'=\bigcup_{k=0}^\infty U_k$. It follows immediately from (1) that $U\subset U'\subset V$. From condition (1), (2) and (3) it follows that
$U_{k+r}\cup V_k = U_k\cup V_k$ for all $k$ and $r>0$, so in particular $U_k\subset U'\subset U_k\cup V_k$ for all $k$. It follows that
\[\del U' \subset \quer{U_k\cup V_k}\setminus U_k \subset \quer{V}_k.\]
Since $V_k$ is $(F_k,d)$-disjoint, we have that $\del U'$ is $(F_k,d)$-disjoint for all $k$, hence it is $(G,d)$-disjoint.
\end{proof}


\section{The controlled marker property}
\noindent
The aim of this section is to use the bounded topological small boundary property to obtain a generalization of Gutman's marker property (see \ref{Gutman}) for free countable group actions and free $\IZ^m$-actions in particular. Very similarly to \ref{dimrok endlich Z}, it will follow that free $\IZ^m$-actions have finite Rokhlin dimension.
First we have  to introduce the notion of markers, the marker property, and, in particular, the controlled marker property.

\defi Let $(X,\alpha,G)$ be a topological dynamical system and $F\fin G$ a finite subset. We call an open set $Z\subset X$ an $F$-marker, if
\begin{itemize}
\item The family of sets $\set{\alpha^g(\quer{Z}) | g\in F}$ is pairwise disjoint.\\ (Or in the notation of the previous section, $\quer{Z}$ is $(F,1)$-disjoint.)
\item $\displaystyle X=\bigcup_{g\in G} \alpha^g(Z)$.
\end{itemize}
We say that $\alpha$ has the marker property if there exist $F$-markers for all $F\fin G$.

\defi Let $(X,\alpha,\IZ^m)$ be a topological dynamical system, $F\fin\IZ^m$ a finite subset and $L\in\IN$ a natural number. We call $Z$ an $L$-controlled $F$ marker, if $Z$ is an $F$-marker such that the second condition can be strengthened to
\[X=\bigcup_{l=1}^L \bigcup_{v\in F} \alpha^{v_l+v}(Z)\quad\text{for certain}\; v_1,\dots,v_L\in \IZ^m.\]
We say that $\alpha$ has the $L$-controlled marker property, if there exist $L$-controlled $B_n$-markers for all $n$.\vspace{3mm}

The marker property has been shown for all aperiodic homeomorphisms (i.e. free $\IZ$-actions) on finite dimensional spaces by Gutman in \cite{Gutman}.
Careful reading of that proof, however, yields the controlled marker property in the form of \ref{Gutman}. It is important to note that although the marker property is trivial if the action is assumed to be minimal, Gutman's proof gives a uniform bound (in relation to $F$) of how many copies one needs to cover the space with an $F$-marker, which is something new even in the minimal case.
We would like to build on his ideas in the case $G=\IZ$ to construct a proof for the general case of countable group actions.

{\lemma \label{key lemma} Let $X$ be a compact metric space, $G$ a countably infinite group, $d\in\IN$ a natural number and $\alpha: G\curvearrowright X$ a free continuous action that satisfies \hyperref[TSBPd]{\em(TSBP $\leq d$)}. 

Let $F\fin G$ a finite subset and let $g_1,\dots,g_d\in G$ be group elements with the property that the sets
\[F^{-1}F\;,\; g_1 F^{-1}F\;,\;\dots\;,\; g_d F^{-1}F\]
are pairwise disjoint. Using the notation $g_0=e$, set $ M=\bigcup_{l=0}^d g_l F^{-1}F$.

Let $U,V\subset X$ be open sets such that
\begin{itemize}
\item The family of sets $\set{\alpha^g(\quer{U}) \;|\; g\in F}$ is pairwise disjoint.
\item The family of sets $\set{\alpha^g(\quer{V}) \;|\; g\in M^{-1}}$ is pairwise disjoint.
\end{itemize}
Then there exists an open set $W\subset X$ such that $U\subset W, V\subset\bigcup_{g\in M} \alpha^g(W)$ and the family of sets $\set{\alpha^g(\quer{W}) \;|\; g\in F}$ is pairwise disjoint.}
\begin{proof}
First observe that there is $\eps>0$ such that $\quer{B}_\eps(U)$ is $(F,1)$-disjoint. Applying \hyperref[TSBPd]{(TSBP $\leq d$)}, we can enlarge $U$ to a new open set within $B_\eps(U)$ such that $\del U$ is $(M,d)$-disjoint. Showing the lemma for this larger set gives certainly no loss of generality, so we may just assume that $\del U$ was $(M,d)$-disjoint to begin with. 

Set $R=\quer{V}\setminus\bigcup_{g\in M} \alpha^g(U)$. Observe that $R$ is closed and $(M^{-1},1)$-disjoint, so choose $\rho>0$ such that $\quer{B}_\rho(R)$ is $(M^{-1},1)$-disjoint as well.
\phantomsection \label{stern2} We now claim that there exists a $\delta>0$ such that
\[(\star)\qquad\qquad |\set{g\in M|\; \alpha^g(\quer{U})\cap \quer{B}_\delta(x)\neq\emptyset}| \leq d\quad\text{for all}\; x\in R.\]
Assume that this is not true. Let $x_n\in R$ be elements with $\delta_n>0$ such that $\delta_n\to 0$ and
\[|\set{g\in M|\; \alpha^g(\quer{U})\cap \quer{B}_{\delta_n}(x_n)\neq\emptyset}|\geq d+1\quad\text{for all}\; n.\]
By compactness, we can assume that $x_n$ converges to some $x\in R$ by passing to a subsequence. Moreover, since $M$ has only finitely many subsets, we can also assume (again by passing to a subsequence if necessary) that there are distinct $\gamma(0),\dots,\gamma(d)\in M$ such that $\alpha^{\gamma(l)}(\quer{U})\cap \quer{B}_{\delta_n}(x_n)\neq\emptyset$ for all $n$ and all $l=0,\dots,d$. But then $\delta_n\to 0$ implies
\[x\in R\cap \bigcap_{l=0}^d \alpha^{\gamma(l)}(\quer{U}) \subset \bigcap_{l=0}^d \alpha^{\gamma(l)}(\del U) = \emptyset.\]
So this gives a contradiction to $\del U$ being $(M,d)$-disjoint.
So we may choose a number $\delta\leq\rho$ satisfying $(\star)$. Moreover, choose some finite covering
\[R\subset\bigcup_{i=1}^s B_\delta(z_i)\quad\text{for some}\; z_1,\dots,z_s\in R.\]
Note that the right side is $(M^{-1},1)$-disjoint by our choice of $\rho$. Now observe that \hyperref[stern2]{$(\star)$} and the fact that the sets $\set{g_l F^{-1}F \;|\; l=0,\dots,d}$ are pairwise disjoint, enables us to define a map $c: \set{1,\dots,s}\to\set{0,\dots,d}$ such that
\[\alpha^g(\quer{U})\cap \quer{B}_\delta(z_i) = \emptyset\quad\text{for all}\; g\in g_{c(i)} F^{-1}F .\]
Finally, set
\[W=U\cup\bigcup_{i=1}^s \alpha^{g_{c(i)}^{-1}}(B_\delta(z_i)).\]
Obviously we have $U\subset W$. Moreover, we have
\[\begin{array}{ccc}
V &\subset& \displaystyle \bigcup_{g\in M} \alpha^g(U) \cup R \\\\
&\subset& \displaystyle  \bigcup_{g\in M} \alpha^g(U) \cup \bigcup_{i=1}^s \underbrace{B_\delta(z_i)}_{=\alpha^{g_{c(i)}}(\alpha^{g_{c(i)}^{-1}}(B_\delta(z_i)))} \\\\
&\subset& \displaystyle  \bigcup_{g\in M} \alpha^g(U) \cup \bigcup_{i=1}^s \alpha^{g_{c(i)}}(W) \quad \subset \bigcup_{g\in M} \alpha^g(W)
\end{array}\]
At last we have to show that $\quer{W}$ is $(F,1)$-disjoint. Suppose that $\alpha^a(\quer{W})\cap\alpha^b(\quer{W})\neq\emptyset$ for some $a\neq b$ in $F$. That is, there exist $x,y\in\quer{W}$ such that $\alpha^a(x)=\alpha^b(y)$. Let us go through all the possible cases:
\begin{itemize}
\item $x,y\in \quer{U}$ is obviously impossible.
\item $x\in\alpha^{g_{c(i_1)}^{-1}}(\quer{B}_\delta(z_{i_1}))$ and $y\in\alpha^{g_{c(i_2)}^{-1}}(\quer{B}_\delta(z_{i_2}))$ for some $1\leq i_1, i_2\leq s$. It follows that
\[\alpha^a(x)=\alpha^b(y)\in \alpha^{ag_{c(i_1)}^{-1}}(\quer{B}_\delta(z_{i_1}))\cap\alpha^{bg_{c(i_2)}^{-1}}(\quer{B}_\delta(z_{i_2})),\]
so
\[\begin{array}{ccl}
\emptyset &\neq& \alpha^{b^{-1} a g_{c(i_1)}^{-1}}(\quer{B}_\delta(z_{i_1}))\cap \alpha^{g_{c(2)}^{-1}}(\quer{B}_\delta(z_{i_2}))\\\\
&\subset& \alpha^{b^{-1}a g_{c(i_1)}^{-1}}(\quer{B}_\rho(R)))\cap \alpha^{g_{c(i_2)}^{-1}}(\quer{B}_\rho(R)).
\end{array}\]
Observe that by $a\neq b$, we have $b^{-1}ag_{c(i_1)}^{-1}\neq g_{c(i_2)}^{-1}$ in $M^{-1}$.
Since $\quer{B}_\rho(R)$ is $(M^{-1},1)$-disjoint, the right side of the above is empty. So this is impossible.
\item $x\in\quer{U}$ and $y\in\alpha^{g_{c(i)}^{-1}}(\quer{B}_\delta(z_i))$ for some $1\leq i\leq s$. Then it follows that
\[\alpha^a(x)=\alpha^b(y)\in\alpha^{a}(\quer{U})\cap\alpha^{bg_{c(i)}^{-1}}(\quer{B}_\delta(z_i))\neq\emptyset.\]
Or equivalently, $\alpha^{g_{c(i)}b^{-1}a}(\quer{U})\cap \quer{B}_\delta(z_i)\neq\emptyset$, a contradiction to the definition of $c(i)$.
\end{itemize}
So the sets in the family $\set{\alpha^g(\quer{W}) | \; g\in F}$ are indeed pairwise disjoint.
\end{proof}

{\prop \label{marker lemma} Let $G$ be a countable group, $X$ a compact metric space, $d\in\IN$ a natural number and $\alpha: G\curvearrowright X$ be a free continuous group action satisfying \hyperref[TSBPd]{\em(TSBP $\leq d$)}. Let $F,g_1,\dots,g_d,M$ be as in \ref{key lemma}. Then $X$ admits an $F$-marker $Z$ with the property $X=\bigcup_{g\in M} \alpha^g(Z)$.}
\begin{proof}
For all $x\in X$, choose a neighbourhood $U_x$ such that the family of sets $\set{\alpha^g(\quer{U}_x) |\; g\in M^{-1}}$ is pairwise disjoint. Note that this is possible because the action is free.

Choose a finite subcovering $X=\bigcup_{i=0}^s U_i$. Apply \ref{key lemma} (with respect to $U=U_0, V=U_1$) to find an open set $W_1$ such that $U_0\subset W_1, U_1\subset\bigcup_{g\in M} \alpha^g(W)$ and such that $\quer{W}_1$ is $(F,1)$-disjoint. Clearly we have $U_0\cup U_1\subset\bigcup_{g\in M} \alpha^g(W_1)$.

Now carry on inductively. If $W_k$ is already defined, apply \ref{key lemma} (with respect to $U=W_k, V=U_{k+1}$) to find $W_{k+1}$ such that $W_k\subset W_{k+1}$ and $U_{k+1}\subset \bigcup_{g\in M} \alpha^g(W_{k+1})$ and such that $\quer{W}_{k+1}$ is $(F,1)$-disjoint. Note also that if $W_k$ had the property that 
\[U_0\cup\dots\cup U_k\subset\bigcup_{g\in M}\alpha^g(W_k),\]
then it follows that
\[\begin{array}{ccl} U_0\cup\dots\cup U_k\cup U_{k+1} &\subset& \bigcup_{g\in M} \alpha^g(W_k) \cup U_{k+1} \\\\
&\subset& \bigcup_{g\in M} \alpha^g(W_k) \cup \bigcup_{g\in M} \alpha^g(W_{k+1}) \\\\
&=& \bigcup_{g\in M} \alpha^g(W_{k+1}).
\end{array}\]
So set $Z=W_s$. $\quer{Z}$ is $(F,1)$-disjoint by construction, and indeed an $F$-marker with 
\[ X=U_0\cup\dots\cup U_s\subset\bigcup_{g\in M} \alpha^g(Z).\]
\end{proof}

\begin{samepage}\rem Let $m\in\IN$. For all $n$ and $a\in\set{0,1}^m$, defining
\[w_a = \Bigl( \delta_{1,a_j}\cdot n \Bigl)_{j=1,\dots,m}\quad\in\quad\IZ^m\]
yields $2^m$ distinct elements with the property that $\displaystyle B_{2n} = \bigcup_a (w_a+B_n)$.\end{samepage}

{\theorem \label{comarkers} Let $X$ be a compact metric space, $d\in\IN$ a natural number and let $\alpha: \IZ^m \curvearrowright X$ be a free continuous action satisfying \hyperref[TSBPd]{\em(TSBP $\leq d$)}.
Then $(X,\alpha,\IZ^m)$ has the $2^m(d+1)$-controlled marker property.}
\begin{proof}
For $n\in\IN$, choose $v_1,\dots,v_d\in\IZ^m$ such that (note $v_0:=0$)
\[B_n-B_n\;,\; v_1+(B_n-B_n)\;,\;\dots\;,\; v_d+(B_n-B_n)\]
are pairwise disjoint. Define $M=\bigcup_{l=0}^d (v_l+(B_n-B_n))$. Apply \ref{marker lemma} to find a $B_n$-marker $Z$ such that $X=\bigcup_{v\in M} \alpha^v(Z)$. Observing that $B_{2n}$ is a translate of $B_n-B_n$, use the previous remark to choose $w_1,\dots,w_{2^m}\in\IZ^m$ so that $B_n-B_n=\bigcup_j (w_j+B_n)$. It follows that
\[M=\bigcup_{l=0}^d (v_l+(B_n-B_n))=\bigcup_{l=0}^d \bigcup_{j=1}^{2^m}((v_l+w_j)+B_n),\]
so
\[X=\bigcup_{v\in M}\alpha^v(Z)= \bigcup_{l=0}^d\bigcup_{j=1}^{2^m} \bigcup_{v\in B_n} \alpha^{(v_l+w_j)+v}(Z).\]
So we see that $Z$ is a $2^m(d+1)$-controlled $B_n$-marker.
\end{proof}


\section{Conclusion and further remarks}
\noindent
Finally we can combine the results from the previous sections to get our main result.

{\cor \label{dimrok endlich} Let $X$ be a compact metric space, $d\in\IN$ a natural number and let $\alpha: \IZ^m \curvearrowright X$ be a free continuous action satisfying \hyperref[TSBPd]{\em(TSBP $\leq d$)}. Then
\[\dimrok(\alpha)\leq 2^m(d+1)-1.\]}
\begin{proof}
Let $n\in\IN$ be given. By \ref{comarkers}, we may apply the $2^m(d+1)$-controlled marker property and choose a $B_n$-marker $Z$ such that there exist 
\[\set{v_l \;|\; l=0,\dots,2^m(d+1)-1}\subset\IZ^m\]
with
\[X=\bigcup_{l=0}^{2^m(d+1)-1} \bigcup_{v\in B_n} \alpha^{v_l+v}(Z).\]
For $l=0,\dots,2^m(d+1)-1$ and $v\in B_n$, define $U_v^{(l)} = \alpha^{v_l+v}(Z).$ Then
\[\CR = \set{ U_v^{(l)}\;|\; l=0,\dots,2^m(d+1)-1,\; v\in B_n}\]
obviously forms a Rokhlin cover with the desired properties.
\end{proof}

{\cor \label{dimrok endlich2} Let $X$ be a compact metric space of finite covering dimension and let $\alpha: \IZ^m \curvearrowright X$ be a free continuous action. Then
\[\dimrok(\alpha)\leq 2^m(\dim(X)+1)-1.\]}
\begin{proof}
This follows directly from \ref{TSBP} and \ref{dimrok endlich}.
\end{proof}

{\theorem \label{main result} Let $X$ be a compact metric space of finite covering dimension and let $\alpha: \IZ^m\curvearrowright X$ be a free continuous group action. Then the induced \cstar-algebraic action $\bar{\alpha}$ on $\CC(X)$ has finite Rokhlin dimension, and the transformation group \cstar-algebra $\CC(X)\rtimes_{\bar{\alpha}}\IZ^m$ has finite nuclear dimension. More specifically, we have
\[\dimrokcyc(\bar{\alpha})\leq 2^{2m}(\dim(X)+1)-1\]
and
\[\dimnuc(\CC(X)\rtimes_{\bar{\alpha}}\IZ^m)\leq 2^{3m}(\dim(X)+1)^2-1.\]}
\begin{proof}
Combining \ref{dimrok endlich2} and \ref{cyc und notcyc} yields
\[\dimrokcyc(\bar{\alpha}) \leq 2^{2m}(\dim(X)+1)-1.\]
Remembering that $\dimnuc(\CC(X))=\dim(X)$, we apply \ref{Hauptsatz} to obtain
\[\dimnuc(\CC(X)\rtimes_{\bar{\alpha}}\IZ^m)\leq 2^{3m}(\dim(X)+1)^2-1.\]
\end{proof}

Now let us adress possible further generalizations and open questions.

\rem The results of this paper apply more generally to free $\IZ^m$-actions on \emph{locally} compact metric spaces. Since the author is not aware of immediate applications for open problems in the non-compact case, that case was dropped to make proofs shorter. However, this may change in the future, so we sketch how one would have to modify the approach to handle the non-compact case.
\begin{itemize}
\item One can define an analogous notion of \cstar-Rokhlin dimension for actions on non-unital \cstar-algebras. It is defined just as in \ref{dimrok}, with the exception that condition (1) is replaced by
\[\hspace{12mm}(1')\quad\|a-a\Bigl(\sum_{l=0}^d\sum_{v\in B_n} f_v^{(l)}\Bigl)\|\leq\eps\quad\text{for all}\; l=0,\dots,d,\; v\in B_n,\; a\in F.\]
In other words, instead of the Rokhlin elements almost summing up to $\eins$, the sum behaves like an approximate unit on elements of $F$.
\item \ref{Hauptsatz} generalizes to actions on non-unital \cstar-algebras. In the proof, condition (1) is used only once, and it is obvious that (1') is actually sufficient in that step.
\item Similarly, one can define an analogous notion of Rokhlin dimension for topological actions on locally compact metric spaces. It is defined just as in \ref{dimrok top}, with the exception that we only require $\CR$ to cover an arbitrarily given compact subset $K\subset X$.
\item \ref{schwacher dimrok} and \ref{cyc und notcyc} can then be generalized in a straightforward manner.
\item Definitions \ref{disjointness}, \ref{SBP} and \ref{genposition} also make sense for actions on non-compact spaces. Note that in \ref{SBP}, the set denoted by $A$ should be required to be compact rather than just closed.
\item Every single statement of section 3 is true in the locally compact case with the same proofs, if we restrict them to statements about bounded open sets. Therefore \ref{TSBP} generalizes to countable group actions on locally compact metric spaces having at most zero dimensional fixed point sets.
\item For group actions on locally compact spaces, we can modify the definition of markers. For $F\fin G$ and a compact subset $K\subset X$, an $(F,K)$-marker is a bounded open set $Z\subset X$ such that
	\begin{itemize}
	\item $\alpha^g(\quer{Z})\cap\alpha^h(\quer{Z})=\emptyset$ for all $g\neq h$ in $F$.
	\item $K\subset\bigcup_{g\in G} \alpha^g(Z).$
	\end{itemize} 
Similarly we can define controlled $(F,K)$-markers in the case of $\IZ^m$-actions. We say that an action has the $L$-controlled marker property if for all $n$ and all compact $K\subset X$, there exists an $L$-controlled $(B_n,K)$-marker.
\item \ref{key lemma}, \ref{marker lemma} and \ref{comarkers} are generalized to the non-compact case in a straightforward manner, using (TSBP $\leq d$) for free actions on locally compact metric spaces. Hence all free $\IZ^m$-actions with (TSBP $\leq d$) have the $2^m(d+1)$-controlled marker property.
\item Combining all these generalized statements, we get the corresponding statements of \ref{dimrok endlich}, \ref{dimrok endlich2} and \ref{main result} for the non-compact case. 
\end{itemize}

\rem Let $G$ be a countable and locally finite group, i.e. $G=\bigcup_n H_n$ for an increasing sequence of finite subgroups $H_n$. Then the results of this paper also apply to free $(\IZ^m\times G)$-actions on (locally) compact metric spaces. Since the generalization is very straightforward, we do not give a sketch on how to modify all the important steps like above. Instead, we remark that the set $B_n$ has to be replaced by $B_n\times H_n$ in order to get the right definitions for $(\IZ^m\times H)$-actions. The proofs are then almost identical.

\question What is the right notion of Rokhlin dimension for actions of a larger class of countable, non-abelian, discrete, amenable groups?

Keep in mind that finite Rokhlin dimension should be regarded as a topological version of the measure theoretic Rokhlin lemma. In view of existing results on groups that can satisfy such a Rokhlin lemma (see \cite{OrnWeiss} in particular), can we define Rokhlin dimension for actions of monotilable groups?

\question Suppose that one can successfully generalize Rokhlin dimension to actions of a larger class of countable, discrete, amenable groups. For which groups is it automatic that free actions on finite dimensional spaces have finite Rokhlin dimension? Lemma \ref{key lemma} seems to suggest that this might be the case for finitely generated groups whose unit balls in some length metric have some form of bounded growth condition similar to $\IZ^m$. Is polynomial growth sufficient?

\end{document}